# Riemann's Zeta Function.
# Numerical Evaluation via its Alternating Relative η(s).

Renaat Van Malderen

September 2011

*Abstract*

The paper describes a method for calculating values of Riemann's Zeta function within the critical strip 0< σ <1 and on its boundary. The approach is based on the "Alternating Zeta function" η(s). The actual Riemann Zeta function ζ(s) is easily obtained from η(s). The obtained accuracy, within certain limits of the described method is discussed. A decent scientific calculator suffices to carry out the involved computations.

*1. Introduction.*

The present paper is not concerned with methods for finding non-trivial zeros of ζ(s). It deals with simply calculating ζ(s) values for s=σ+it in the range 0< σ <1 . The method is also applicable to "almost all" points on the verticals s=it and s=1+it with the exception of a discrete set of points as will be explained in section 5. H.M. Edwards in his wonderful book "Riemann's Zeta Function" describes the "Euler-Maclaurin Summation" approach [1, p.98] for computing ζ(s).

A similar approach also starting from the basic series

$$\zeta(s) = \sum_{n=1}^{\infty} n^{-s} \text{ for } \sigma > 1 \qquad (1)$$

deals with the simple pole of ζ(s) at s=1 by splitting the series (1) into an integral and a remaining difference. The resulting expression allows to analytically continue ζ(s) to the entire complex plane. One example of such an approach expresses (1) as:

$$\zeta(s) = 1 + \int_2^{\infty} \frac{dx}{\left(x - \frac{1}{2}\right)^s} - \sum_{n=2}^{\infty} \int_0^1 \left[ \frac{1}{\left(n - \frac{1}{2} + \lambda\right)^s} - \frac{1}{n^s} \right] d\lambda \qquad (2)$$

which after some algebra yields:

$$\zeta(s) = 1 + \frac{1}{s-1}\left(\frac{2}{3}\right)^{s-1} - \sum_{k=1}^{\infty} C_k(s)[\zeta(s+2k) - 1] \qquad (3)$$

with coefficients:

$$C_k(s) = \frac{\prod_{r=0}^{2k-1}(s+r)}{(2k+1)!\, 2^{2k}} \quad (4)$$

(3) is absolutely convergent and valid for all complex s (except for s=1 of course) whereas in the Euler-Maclaurin approach the involved series ultimately "grows without bound" [1, p. 105].
(3) however requires the knowledge of the terms (ζ(s+2k)-1).

For s values with $0 \leq \sigma \leq 1$ obtaining the first few terms via series (1) requires to compute a significant number of terms in (1) due to slow convergence. Also, for increasing t values, more and more $C_k(s)$ need to be taken into account before the terms in (3) start to taper off. As a side issue, for σ +2k values ≤ 1, series (1) does not converge anymore at all and (3) has to be used repeatedly, moving stepwise to the left.

*2. The Alternating Zeta Function η(s).*

η(s) is defined as:

$$\eta(s) = \sum_{n=1}^{\infty} (-1)^{n+1} n^{-s} \quad (5)$$

with the complex variable s= σ +it. (5) is convergent for σ >0. As well known, the actual Riemann Zeta function ζ(s) may be obtained from (5) through the relation:

$$\zeta(s) = \frac{\eta(s).}{1 - 2^{1-s}} \quad (6)$$

ζ(s) is analytic in the entire s plane except for the simple pole at s=1. At that point η(1)=ln2 as can be seen from (5). All the zeros of ζ(s) are also zeros of η(s), both the trivial ones as well as those occurring in the critical strip 0< σ <1.

In addition η(s) has simple zeros for:

$$s = 1 \pm \frac{2\pi n i}{\ln 2} \text{ with } n = 1, 2, 3, \dots \quad (7)$$

This is obvious from (6) and the fact that ζ(s) has no zeros on the line σ =1.
η(s) may also be represented as an integral valid for σ >0 [2, p.267]:

$$\eta(s) = \frac{1}{\Gamma(s)} \int_0^\infty \frac{\lambda^{s-1} d\lambda}{1+e^\lambda} = \frac{1}{\Gamma(s)} \int_0^\infty \sum_{n=1}^\infty \lambda^{s-1} e^{-n\lambda}(-1)^{n+1} d\lambda$$

$$\eta(s) = \frac{1}{\Gamma(s)} \sum_{n=1}^{2m} \int_0^\infty \lambda^{s-1} e^{-n\lambda}(-1)^{n+1} d\lambda + \frac{1}{\Gamma(s)} \sum_{n=2m+1}^\infty \int_0^\infty \lambda^{s-1} e^{-n\lambda}(-1)^{n+1} d\lambda \quad (8)$$

Through substitution $n\lambda = x$ in the first integral and since $\int_0^\infty x^{s-1} e^{-x} dx = \Gamma(s)$:

$$\eta(s) = \sum_{n=1}^{2m} (-1)^{n+1} n^{-s} + \frac{1}{\Gamma(s)} \int_0^\infty \frac{\lambda^{s-1} e^{-(2m+1)\lambda} d\lambda}{1+e^{-\lambda}}$$

$$\eta(s) = \sum_{n=1}^{2m} (-1)^{n+1} n^{-s} + I_1(s, m) \quad (9)$$

Via substitution $(2m+1)\lambda = y$ and putting $\beta = 1/(2m+1)$ we obtain:

$$I_1(s, \beta) = \frac{(2m+1)^{-s}}{\Gamma(s)} \int_0^\infty \frac{y^{s-1} e^{-y} d\lambda}{1+e^{-\beta y}} \quad (10)$$

$\frac{1}{1+e^{-\beta y}}$ may be expanded in a power series within the convergence circle $|\beta y| < \pi$:

$$\frac{1}{1+e^{-\beta y}} = \frac{1}{2}\left[1 + \text{Th}\left(\frac{\beta y}{2}\right)\right]$$

As well known:

$$\text{Th}\left(\frac{\beta y}{2}\right) = \sum_{r=1}^\infty \frac{2^{2r}(2^{2r}-1)B_{2r}(\beta y)^{2r-1}}{(2r)! \, 2^{2r-1}}$$

where $B_{2r}$ is the Bernoulli number of index $2r$ and $\left|\frac{\beta y}{2}\right| < \frac{\pi}{2}$.

$$\frac{1}{1+e^{-\beta y}} = \frac{1}{2}\left[1 + \sum_{r=1}^\infty \frac{2(2^{2r}-1)B_{2r}(\beta y)^{2r-1}}{(2r)!}\right] \text{ for } |\beta y| < \pi$$

$$\frac{1}{1+e^{-\beta y}} = P_0 + \sum_{r=1}^\infty P_{2r-1}(\beta y)^{2r-1} \quad (11)$$

with $P_0 = \frac{1}{2}$ and $P_{2r-1} = \frac{(2^{2r}-1)B_{2r}}{(2r)!}$ for r=1, 2, 3, ... (12)

## 2. Coefficient Values

Table 1 gives coefficients for the series (11).
Using:
$\frac{B_{2r}}{(2r)!} \cong \frac{2(-1)^{r+1}}{(1-2^{1-2r})(2\pi)^{2r}}$ we obtain for large r:

$P_{2r-1} \cong \frac{2(-1)^{r+1}}{\pi^{2r}}$ (13)

Table 1

|       |                     |                              |
|-------|---------------------|------------------------------|
| $P_0$ | 1/2                 | 0.5                          |
| $P_1$ | 1/4                 | 0.25                         |
| $P_3$ | -1/48               | -2.083333.... x $10^{-2}$    |
| $P_5$ | 1/480               | 2.083333.... x $10^{-3}$     |
| $P_7$ | -17/80640           | -2.10813492 x $10^{-4}$      |
| $P_9$ | 31/1451520          | 2.1356922 x $10^{-5}$        |
| $P_{11}$ | -691/319334400   | -2.163876 x $10^{-6}$        |
| $P_{13}$ | 5461/24908083200 | 2.19246 x $10^{-7}$          |

## 3. Approximations

We now discuss the approximations we will make in order to obtain numerical values for $I_1(s,\beta)$ given in (10). For convenience, wherever appropriate, instead of $\beta$ we will use $b = \frac{1}{\beta} = 2m + 1$.

In addition we will stay well within the convergence circle (11) by imposing $\beta y \leq 1$.

With this in mind $I_1(s, \beta)$ equals:

$$I_1(s,\beta) = \frac{(2m+1)^{-s}}{\Gamma(s)}\left[\int_0^b y^{s-1}e^{-y}dy\left(P_0 + \sum_{r=1}^{\infty}P_{2r-1}(\beta y)^{2r-1}\right) + \int_b^{\infty}\frac{y^{s-1}e^{-y}dy}{1+e^{-\beta y}}\right] \quad (14)$$

Important restriction: From now on our considerations will be limited to $s = \sigma + it$ values with $\frac{1}{2} \leq \sigma \leq 1$ (15).

Riemann's functional equation [3, p. 438] allows to obtain from $\zeta(\frac{1}{2} + \delta + it)$ first $\zeta(\frac{1}{2} - \delta - it)$, and then $\zeta\left(\frac{1}{2} - \delta + it\right) = \overline{\zeta(\frac{1}{2} - \delta - it)}$.

In this way the range $0 \leq \sigma \leq 1$ is covered.

To proceed we will make three approximations regarding the integral $I_1(s,\beta)$:

a) We limit the number of terms in the power series for $\beta y$ to a maximum value $r=r_m$.

b) By selecting large enough values for $b = \frac{1}{\beta}$ (i.e. for m) we intend to ignore $\int_b^{\infty}$ in (14).

c) The integrals $\int_0^b$ in (14) are replaced by $\int_0^{\infty}$.

The error quantification and impact on the calculations are dealt with below. Due to the limitation (15) on $\sigma$ we can easily give an upper limit for the absolute value of the last integral in (14):

$$\left|\int_b^{\infty}\frac{y^{s-1}e^{-y}dy}{1+e^{-\beta y}}\right| < \int_b^{\infty}e^{-y}dy = e^{-b}$$

or

$$\int_b^{\infty}\frac{y^{s-1}e^{-y}dy}{1+e^{-\beta y}} = \vartheta_1(s)e^{-b} \quad (16)$$

with $\vartheta_1(s)$ generally complex and $|\vartheta_1(s)| < 1$.

For m = 40 as used in some of the numerical examples, $e^{-(2m+1)} = e^{-81} < 7 \times 10^{-36}$.

By limiting the expression of $(1 + e^{-\beta y})^{-1}$ to the term $(\beta y)^{13}$, i.e. $r_m = 7$:

$$\frac{1}{1+e^{-\beta y}} \cong P_0 + \sum_{r=1}^{7} P_{2r-1}(\beta y)^{2r-1} \quad (17)$$

we achieve an accuracy in $(1 + e^{-\beta y})^{-1}$ of :
  error $< 10^{-8}$ for $\beta y \leq 0.95$
  error $< 2 \times 10^{-8}$ for $\beta y \leq 1$

This may be verified by actual calculation and the resulting error is independent of $\beta$ as long as $\beta y \leq 1$.
Taking the above into account:

$$I_1(s,\beta) = \frac{(2m+1)^{-s}}{\Gamma(s)} \left[ P_0 \left( \Gamma(s) - \int_b^\infty y^{s-1} e^{-y} dy \right) \right.$$

$$+ \sum_{r=1}^{7} P_{2r-1} \beta^{2r-1} \left( \Gamma(s+2r-1) - \int_b^\infty y^{s-1+(2r-1)} e^{-y} dy \right)$$

$$\left. + \theta(s) e^{-b} \right] \quad (18)$$

As the next step we want to drop the integrals $\int_b^\infty y^{s-1+n} e^{-y} dy$ with n=0, 1, 3, 5, 7, 9, 11, 13 in (18).
This requires an estimate of the resulting error. So we need to compare:

$|\Gamma(s+n)|$ against $\left| \int_b^\infty y^{s-1+n} e^{-y} dy \right|$. (19)

The variables involved are $s = \sigma + it$ (with $\frac{1}{2} \leq \sigma \leq 1$) and b= 2m+1.
$|\Gamma(s+n)| = |\Gamma(\sigma + n + it)|$ grows with n but decreases with t. It is not easy to figure out $\left| \int_b^\infty y^{\sigma-1+n+it} e^{-y} dy \right|$ for $t \neq 0$.
However $\left| \int_b^\infty y^{\sigma-1+n+it} e^{-y} dy \right| \leq \int_b^\infty y^{\sigma-1+n} e^{-y} dy$ (20).

Let's first consider the case $\sigma = 1$. (20) then results in a solvable integral:

$$\left| \int_b^\infty y^{n+it} e^{-y} dy \right| \leq \int_b^\infty y^n e^{-y} dy = E(n,b) = e^{-b} n! \sum_{k=0}^{n} \frac{b^k}{k!} \quad (21)$$

The role of the parameter t:

Turning now to:

$$|\Gamma(s+n)| = \left|\int_0^\infty y^{\sigma-1+n+it}e^{-y}dy\right| = \left|\int_0^\infty y^{\sigma-1+n}\,e^{-y}e^{it\ln y}\,dy\right| \quad (22)$$

Consider the last integral in (22). Its magnitude is determined by two aspects: Due to the oscillatory nature of $e^{it\ln y}$ one would expect the main contribution to (22) to come from the y-range where the amplitude and the derivative of the curve $y^{\sigma-1+n}\,e^{-y}$ is most significant.
However, the higher the "frequency" of $e^{it\ln y}$ for a given y, the lower the contribution to (22) at that point. For y→ 0, this frequency goes to ∞, and for increasing y, this frequency keeps slowing down. The overall result will be to shift the y-area where the bulk of (22) will be picked up, to the right with increasing t. Unless we select b large enough, a significant amount of (22) will, for larger and larger t, eventually be picked up in the range y > b.

Considering first again σ = 1.

$$\left|\int_0^\infty y^n e^{-y}e^{it\ln y}dy\right| = |\Gamma(n+1+it)|$$

Using Stirling's formula for x+it:

$$|\Gamma(x+it)| \cong \sqrt{2\pi}\,e^{-x}[x^2+t^2]^{\left(\frac{x}{2}-\frac{1}{4}\right)}e^{-t\,\text{arctg}\left(\frac{t}{x}\right)} \quad (23)$$

we can calculate $|\Gamma(n+1+it)|$ with excellent accuracy and compare it with the worst case error (t=0) as given by (21).

As pointed out above $|\Gamma(n+1+it)|$ decreases with increasing t. As long as E(n,b) is very small compared to $|\Gamma(n+1+it)|$, replacing $\int_0^b$ by $\int_0^\infty$ in (14) is acceptable. For increasing t the occurring error will keep growing until eventually E(n,b) would be larger than $|\Gamma(n+1+it)|$ which of course makes no sense. Actual calculations and comparisons between $|\Gamma(n+1+it)|$ and E(n,b) show specifically that:
For b= 81, n ≤ 13, t ≤ 40: E(n,b) < $10^{-5}$ $|\Gamma(n+1+it)|$ (24)
For b= 101, n ≤ 13, t ≤ 50: E(n,b) < $10^{-7}$ $|\Gamma(n+1+it)|$ (25)
Again, these are conservative estimates since the reduction of E(n,b) due to t is not taken into account. Calculations for larger t values require correspondingly larger b values (i.e. m in (9) will have to increase).
Finally for the case $\sigma = \frac{1}{2}$, $|\Gamma(n+1/2+it)|$ and the error $\int_b^\infty y^{n-\frac{1}{2}}\,e^{-y}dy$ will both come down with essentially similar results.

## 4. Resulting Expression for $\eta(s)$.

Resulting from the above approximations (18) becomes:

$$I_1(s, \beta) \cong \frac{(2m+1)^{-s}}{\Gamma(s)} \left[ P_0 \Gamma(s) + \sum_{r=1}^{7} P_{2r-1} \beta^{2r-1} \Gamma(s+2r-1) \right]$$

Using the fundamental property of the Gamma function $\Gamma(z+1) = z\, \Gamma(z)$ we can eliminate $\Gamma(s)$ in the numerator and denominator:

$$I_1(s, \beta) \cong (2m+1)^{-s} \left[ P_0 + \sum_{r=1}^{7} P_{2r-1} \beta^{2r-1} \prod_{k=0}^{2(r-1)} (s+k) \right]$$

(9) now becomes:

$$\eta(s) \cong \sum_{n=1}^{2m} (-1)^{n+1} n^{-s} + (2m+1)^{-s} \left[ P_0 + \sum_{r=1}^{7} P_{2r-1} \beta^{2r-1} \prod_{k=0}^{2(r-1)} (s+k) \right] \quad (26)$$

## 5. Numerical results.

This section contains a number of specific examples to demonstrate the accuracy of formula (26). These examples are all dealing with $s = \sigma + it$ values for which $\frac{1}{2} \leq \sigma \leq 1$ and $0 \leq t \leq 50$.

All numerical calculations listed in this section were carried out with a Texas Instruments TI-89 Titanium calculator, displaying 12 digits after the decimal point.

a) $\eta(1) = \ln 2$ as mentioned in section 2.
   ln 2 as given in [4, p.113] up to 13 digits equals
      ln 2 = 0.6931471805599...
   Formula (26) yields:
   For m= 40: $\eta(1)$= 0.693147181 ; Error <1.5 x $10^{-9}$
   For m= 50: $\eta(1)$= 0.69314718056 ; Error <1 x $10^{-11}$
      Value for $\zeta(1/2)$:
   For m= 40: s=0.5, (26) yields $\eta(0.5)$= 0.604898643422
   For m= 50: s=0.5, (26) yields $\eta(0.5)$= 0.604898643422
   $\zeta(1/2) = \frac{\eta(0.5)}{(1-\sqrt{2})} = $ -1.4603545...

b) Zeros of $\eta(s)$ on the line 1+it provide a good method to check the obtained accuracy of (26). Indeed, we know the result ought to be exactly zero. The small difference with zero resulting from (26) shows the actual error. Table 2 gives $\eta(s)$ for $s = 1 + it = 1 + \frac{n\pi i}{\ln 2}$ for values < 50. Even –n values correspond to zeros. The non-zero values (n odd) are halfway between the zeros and were rounded off to six digits after the decimal point. The errors on the zeros of $\eta(s)$ are $< 10^{-12}$.

Table 2: $s = 1 + it = 1 + \frac{n\pi i}{\ln 2}$

| $n$ | $\eta(s)$ | zeros |
|---|---|---|
| 0 | ln 2 | |
| 1 | 1.437551 + 0.249393 i | |
| 2 | -3.9 x $10^{-14}$ + 2.29 x $10^{-14}$ i | X |
| 3 | 0.803791 – 0.442143 i | |
| 4 | -2.1 x $10^{-14}$ – 2.11 x $10^{-14}$ i | X |
| 5 | 2.111898 + 0.441761 i | |
| 6 | -7.31 x $10^{-14}$ + 1.8516 x $10^{-13}$ i | X |
| 7 | 4.893007 – 0.100872 i | |
| 8 | -5.3 x $10^{-14}$ + 8.1 x $10^{-14}$ i | X |
| 9 | 0.973771 – 0.301917 i | |
| 10 | -2.22 x $10^{-13}$ + 1.15 x $10^{-14}$ i | X |
| 11 | 0.855469 + 0.382103 i | |

c) Zeros on the critical axis $\sigma = 1/2$.
Remember, these zeros are the same for $\zeta(s)$ and $\eta(s)$. Two s values were used: the lowest non-trivial zero and the tenth one, just below $s = \frac{1}{2} + 50i$. The s values are rounded-off to six decimal places. Accordingly a value of m=40 was used in (26).
For $s_1 = 0.5 + 14.134725i$ the error compared to zero is:
$\eta(s_1) = 1.6212$ x $10^{-8}$ – 2.6635 x $10^{-7}$ i.
For $s_{10} = 0.5 + 49.773832i$ the error compared to zero is:
$\eta(s_{10}) = 7.98233$ x $10^{-7}$ – 1.427674 x $10^{-6}$ i.
These errors have two causes:
a) the finite accuracy of $s_1$ and $s_{10}$.
b) the error introduced by (26).

d) Riemann's functional equation in terms of η(s) is:

$$\eta(1-s) = \frac{\pi \, \eta(s)(1-2^s)}{\Gamma(1-s)\sin\left(\frac{\pi s}{2}\right)(2\pi)^s(1-2^{1-s})} \quad (27)$$

(27) allows to obtain η(0) from η(1)=ln2:
Plugging – in s=ln2 in (27) requires to consider the limit
for $s = 1 + \varepsilon$ with $\varepsilon \to 0$ of :

$$\lim_{\varepsilon \to 0} \Gamma(-\varepsilon)(1-2^{-\varepsilon}) = -\ln 2$$

(27) then yields: η(0)=1/2. Accordingly ζ(0)=-1/2.
As already hinted at in the introduction, since η(1± 2πni/ln2)=0 we cannot use (27) to obtain η(1-s) since (27) depends on η(s)/(1-2$^{1-s}$) which for s=1± 2πni/ln2 gives a $\frac{0}{0}$ situation. To resolve this we ought to apply l'Hopital's rule, but to figure out $\frac{d\eta(s)}{ds}$ is no straightforward task. From the theory of Dirichlet series ([5, p. 445] or [6, p.236]) we know that differentiating (5) is allowed for σ>0 but the convergence of the new series is even slower than (5).
The better way is to first compute ζ(1± 2πni/ln2) with e.g; formula (3), then ζ(± 2πni/ln2)using Riemann's functional equation in terms of ζ(s) and finally use (6) to obtain η(± 2πni/ln2).

Renaat Van Malderen
Address: Maxlaan 21, B-2640 Mortsel, Belgium
The author can be reached by email: hans.van.malderen@telenet.be